\documentclass[twoside]{amsart}

\usepackage{amsfonts}
\usepackage{amsmath}   %multiline etc. Subsumes amstext,amsbsy,amsopn
\usepackage{amsthm}    %provide proof and exercise env
\usepackage{amscd}     %commutative diagrams

\usepackage{euscript}

\usepackage{amssymb}

\newtheorem{thm}{Theorem}[section]
\newtheorem{lemma}[thm]{Lemma}
\newtheorem{cor}[thm]{Corollary}

\theoremstyle{definition}

\theoremstyle{remark}

\newcommand{\ZZ}{{\mathbb Z}}

\newcommand{\OO}{{\EuScript O}}
\newcommand{\II}{{\EuScript I}}

\newcommand{\pin}{{\bf P}^{n}}

\newcommand{\pro}{{\bf P}^{m+r}}
\newcommand{\frd }{\longrightarrow}
\begin{document}
\title{On quadratic and higher normality of small codimension 
projective varieties}

\author{Chiara Brandigi}

\address{Chiara Brandigi\\
Dipartimento di matamatica U.Dini\\
Viale Morgagni 67/A\\
50134 Firenze Italy}
\email{brandigi@math.unifi.it}

\markboth{CHIARA BRANDIGI}{ON QUADRATIC AND HIGHER NORMALITY OF SMALL
  CODIMENSION VARIETIES}

\begin{abstract}
  Ran proved that smooth codimension $2$ varieties in ${\bf P }^{m+2}$
  are $j$-normal if $(j+1)(3j-1)\le m-1$, in this paper we extend this
  result to small codimension projective varieties. Let $X$ be a $ r$
  codimension subvariety of $\pro$, we prove that if the set
  $\Sigma_{(j+1)}$ of $(j+1)$-secants to $X$ through a generic
  external point is not empty, $ 2(r+1)j\leq m-r $ and $
  (j+1)((r+1)j-1)\leq m-1$ then $X$ is $j$-normal.  If $X$ is given by
  the zero locus of a section of a rank $r$ vector bundle $E$ on
  $\pro$, we prove that $ \textrm{deg }
  \Sigma_{j+1}=\frac{1}{(j+1)!}\prod_{i=0}^{j}c_{r}(E(-i))$. Moreover we get a new
  simple proof of Zak's theorem on linear normality if $m\ge 3r$.
  Finally we prove that if $c_{r}(N(-2))\neq 0$ and $6r\le m-4$ then
  $X$ is $2$-normal.
\end{abstract}  


\maketitle
\tableofcontents

\section{Introduction}
A variety $X\subset {\bf P}^{n}$ is called $j$-normal if the restriction map
$H^{0}({\bf P}^{n},\OO (j))\longrightarrow H^{0}(X,\OO(j))$ is surjective.
Hartshorne's conjecture \cite{H1} implies that smooth varieties $X\subset \pin$ of 
small codimension are $j$-normal.
Peternell, Le Potier , Schneider \cite{P-L-S} and Ein \cite{E} proved indipendently
that smooth codimension $2$ varieties $X\subset \pin$ are $2$-normal if $n\ge 10$. This bound is probably not sharp 
(Hartshorne's conjecture implies  $n\ge 6$) but it is interesting
because it does not depend on the degree of $X$ (for similar bounds 
depending on the degree, see \cite{G-L-P} \cite{K}).Ein's results were extended to higher
normality by Alzati and Ottaviani 
in \cite{A-O}, but the techniques of those  papers seem not to work in codimension 
$\ge 3$ because the Koszul complexes appearing in the proof have greater length and are difficult to control.
On the other hand  Ran in \cite{R2} proved, with different techniques, that smooth codimension $2$ varieties 
$X\subset \pin$ are $j$-normal if $n\ge 3j^{2}+2j+2 $. Ran constructs explicitly, for any $Y\in H^{0}(X,\OO (j))$, a hypersurface $F$ in $\pin$ of degree $j$ as the union of lines which intersects $Y$ with multiplicity $\ge j+1$. This works  because the assumption implies that the locus of $j+1$-secants is not empty.
In our doctoral thesis, we expanded all the details of  Ran's  paper and 
we were able to prove the following theorem which gives bounds 
for $j$-normality also in codimension $r\ge 3$. \\

Denote by $\Sigma_{(j+1)}$ the set of $(j+1)$-secants to $X$ through a (generic) external 
point.

\begin{thm}
\label{T1}
Let  $X$ be a $ r$ codimension subvariety of $\pro$;  if 
\[\Sigma_{(j+1)}\neq \emptyset\]
\[ 2(r+1)j\leq m-r \qquad 
\textrm{and} \qquad (j+1)((r+1)j-1)\leq m-1 \] then:
\[ \rho_{j} : H^{0}(\pro,\OO_{\pro}(j))
\longrightarrow H^{0}(X,\OO_{X}(j)) \]
is  surjective.
\end{thm}
If $r=2$, the numerical assumptions of theorem $1.1$ are exactly as in \cite{R2}, while Ran is able to show that in this bound  if $\Sigma_{j+1}=\emptyset$ then $X$ is a complete intersection. 

Ran himself  pointed out in a remark at the end of the paper that his 
proof could also be extended  to higher codimension.
 When $X$ is the zero locus of a section of a vector bundle,
then the numeric assumption is more explicit.
\begin{thm}
\label{T2}
Let $X$ be a $m$ dimension variety  in $\pro$
given by the zero locus of a section of a rank $r$ vector bundle
$E$ on $\pro$. We have
$$ \textrm{deg } \Sigma_{j+1}=\frac{1}{(j+1)!}\prod_{i=0}^{j}c_{r}(E(-i))$$ 
\end{thm} 
\begin{cor}
With the assumptions of the theorem $1.2$, if
$$ c_{r}(E(-i))\neq 0 \quad \forall i= 1\ldots j$$

\[ 2(r+1)j\leq m-r \qquad 
and \qquad
 (j+1)((r+1)j-1)\leq m-1 \]
then:
\[  \rho_{j} : H^{0}(\pro,\OO_{\pro}(j))
\longrightarrow H^{0}(X,\OO_{X}(j)) 
\]
is surjective.
\end{cor}
In section $4$ we get a new proof of Zak theorem about linear normality with the assumption $n\ge 4r$. In the same range there is still another proof due to Faltings \cite{Fa}.
Moreover in this paper we prove the following result on quadratical 
normality  where the numeric assumption is easier  checked. This is a partial answer to problem $12$ in Schneider list \cite{S}.

%\teo3
\begin{thm}
\label{T3}
Let $X$ be a $m$ dimension variety  in $\pro$. If 
\[ c_{r}(N(-2))\neq 0\quad  and\quad 6r\leq m-4\]
then $X$ is  $2$-normal.
\end{thm}
I thank  G.Ottaviani for the precious suggestions and the useful discussions and  L.G\"ottsche for some ideas used for the proof of theorem $1.2$.

\section{Proof of theorem $1.1$}

\newtheorem{guess}{Esempio}

Consider a branched covering that is a finite surjective morphism 
between two  irreducible and nonsingular algebraic varieties $V$ and $W$  
$f:V\rightarrow W $; let $d$ be the degree of $f$. As we are assuming that $V$ 
and $W$
are non-singular, $f$ is flat and consequently the direct image $f_{*}\OO_{V}$ 
is locally free of rank  $d$ on $W$. The trace $Tr_{V/W}:f_{*}\OO_{V}\rightarrow\OO_{W}$ 
gives rise to a splitting:
$f_{*}\OO_{X}=\OO_{W}\oplus F$, where $F=ker(Tr_{V/W}).$
We shall be concerned with the rank  $d-1$ vector bundle on $W$:
$E=F^{*}.$
E will be termed {\em vector bundle associated with the covering} $f$.
Let
$e_{f}(x)=dim_{{\bf C}}(\OO_{x}X/f^{*}m_{f(V)})$
be the local degree of $f$ in $x$ which counts the number of sheets of covering that
come together at $x$.\\
 
\begin{thm}[Gaffney-Lazarsfeld] Let $V$ and $W$ be varieties of dimension $n$ and $f:V\frd  W$ 
a branched covering of degree $d$; if the vector bundle associated with a branched 
covering is ample, then there exists at least one point $x\in V$ at which
$$ e_{f}(x)\ge min(d, n+1).$$  
\end{thm}
 
{\bf Proof } See \cite{G-L}. Lazarsfeld himself points out that smoothness of $W$ is not essential.\\
                            
Thanks to this  theorem,  we are able to prove the following Lemma:
\begin{lemma}
  Let $X$ be a $r$-codimensional subvariety of ${\bf P}^{n}$, if $ r\cdot k
  \ge n $ and the set of $k$-secant lines to $X$ through an external point $P$ is
  not empty, then there exists at least a $k$-secant  through this point at which the
  $k$ points coincide.
\end{lemma}

{\bf Proof}
 We consider the projection from $P$ of $k$-secants on a generic hyperplane  ${\bf P}^{n-1}$; 
 let $f$ be its restriction to the points of  $X$, and $Y$ the image of $f$.
 $X'=f^{-1}(Y)$,
 $X'$ is the set of points in $X$ lying on a $k$-secant. 
 The dimension $n'$ of  $X'$ and $Y$
is $n-1-k(r-1)$ and
$f:X'\frd  Y$ is a finite covering with degree $k$: by our assumptions 
the degree of the covering is less than or equal to $n'+1.$
If we prove that the vector bundle associated with the  covering is
ample, then we can use the theorem of Gaffney and Lazarsfeld to prove that there exists
a point at which the sheets of covering come together.\ We denote $C$ as the
cone of $k$-secants through an external point $P$; since there are $k$
points of $X'$ for each $k$-secant, we observe that $X'$ is a
divisor of $C$ and since the point $P$ is external to $X$ this divisor
is disjoint from singularities of $C$.  Let $C'$ be the
desingularization of $C$, we have:
$$C'={\bf P}(\OO_{Y}\oplus \OO_{Y}(1)),$$
then $X'$ is
isomorphic to a divisor of $C'$.
$f_{*}\OO_{X'}$ is a vector bundle of rank $k$; we want to prove
that:
$$ f_{*}\OO_{X^{'}}=\OO_{Y}\oplus \OO_{Y}(-1)\oplus\ldots \oplus \OO_{Y}(1-k).$$
$X'$ is the zero locus of a section of $\OO_{{\bf P}(\OO\oplus \OO(1))}(k)$; 
in fact, from \cite{H2} we have: 
$Pic(C') =Pic(Y)\oplus \ZZ H$,
where $H$ is hyperplane section.
$X'$ is a divisor which meets the generic fibre in $k$ points and it is  
disjoint to the infinite section, and so $X$ is linearly equivalent to 
$kH$.
Now we consider the associated exact sequence:
$$0\frd  \OO_{{\bf P}}(-k)\frd  \OO_{{\bf P}}\frd  \OO_{X^{'}}\frd  0$$
Let  $\pi$ be the projection from ${\bf P}(\OO\oplus \OO(1))$ to $Y$; 
applying $\pi _{*}$ to the sequence  we obtain:
$$0\frd  \OO_{Y}\frd  \pi_{*} \OO_{X'}\frd  R^{1}\pi_{*}\OO_{{\bf P}}(-k)\frd  0.$$
Using the exercise $8.4$ of \cite{H2}, (page $253$) 
we prove that:
$$R^{1}\pi_{*}\OO(-k)\simeq \pi_{*}(\OO(k-2))^{*}\otimes \OO_{Y}(-1)$$
and from the same exercise we have:
$$\pi_{*}\OO(k-2)\simeq S^{k-2}(\OO\oplus \OO(1))=\OO\oplus \OO(1)\oplus \ldots \oplus \OO(k-2)$$
then
$$R^{1}\pi_{*}\OO(-k)=\OO(-1)\oplus \OO(-2)\oplus \ldots \OO(-k+1)$$
substituting in the exact sequence we get:
$$\pi_{*}\OO_{X^{'}}=\OO_{Y}\oplus \OO_{Y}(-1)\oplus\OO_{Y}(-2)\oplus \ldots \OO_{Y}(1-k)$$
then
$$\pi_{*}\OO_{X^{'}}=\OO_{Y}\oplus F$$
$$\pi_{*}\OO_{X^{'}}=f_{*}\OO_{X^{'}}$$
where $F$ is a vector bundle whose dual is ample.
We can now use  the Gaffney-Lazarsfeld's theorem to obtain the thesis.

\begin{lemma}
Let $G$ be a generic hypersurface of ${ \bf P}^{n}$ of degree $j$ passing 
through a point $P$, then the variety of lines through $P$ lying  in $G$ is a 
 complete intersection of ${\bf P}^{n-2}$ with dimension $n-j-1$ and degree $j!$.
\end{lemma}
{\bf Proof.}  
We can choose a coordinate system  such that P is the point $(a,0,0,\ldots,0)$.
Let $\pi$ be the hyperplane $x_{0}=0$; for every point $Q$ of $\pi$ we
 consider the line $r$ through $P$ and $Q$ that is $(a(1-t),t
 x_{1},\ldots,t x_{n}).$ G is given by $F(y_{0},\ldots,y_{n})=0$ with
 $F(y_{0},\ldots,y_{n})= by_{0}^{j}+f_{1}(y_{1},\ldots ,y_{n})
 y_{0}^{j-1}+\ldots + f_{j}(y_{1},\ldots ,y_{n}) $ where $f_{i}$ are
 polynomials of degree $i$; since $P\in G$ we have $b=0$.  A line $r$
 lie on G if and only if:
\begin{displaymath}
F(a(1-t),t x_{1},\ldots,t x_{n})=ty_{0}^{j-1}f_{1}(x_{1},\ldots
x_{n})+ \ldots t^{j}f_{j}(x_{1},\ldots ,x_{n})=0 \quad
 \end{displaymath} 
for every $t$, and so we must have: 
$f_{i}(x_{1},\ldots,x_{n})=0\qquad  \forall i=1,\ldots,n.$
Since $G$ is generic and $f_{1}$ is linear, this gives a transversal  intersection  contained in ${\bf P}^{n-2}$.
Finally we get that the variety of lines of $G$ through a point $P\in G$
is a complete intersection of degree $j!$ and dimension $n-1-j$.

Let $X$ be a subvariety of  ${\bf P}^{m+r}$; we denote by $\Sigma_{j} $ 
the cycle of  $j$-secant lines to $X$ through an external point.

{\bf Proof  of theorem 1.1}   Consider a generic element $Y$ of the linear system
$\mid \OO_{X}(j)\mid $.  Since the locus of 
$(j+1)$-secants through a generic point is not empty, then $X$ can not be included
in a hypersurface of degree $j$ and so $H^{0}({\II}_X(j))=0$.  \\
In order to prove the theorem we  just have to find one hypersurface of degree 
 $\leq j$  which contains $Y$.\\

\newcommand{\re}{$R^{k}$}
We define
\re =$\{(y,z)\in Y \times{\bf P}^{m+r}$ : $\exists$ a line L from 
$z\in \pro$ such that $L\cap Y$ has multiplicity $ \ge k$ in $ y$\}.
Let $p$ and $q$ the projections of $R^{k}$  to $ Y$ and to $\pro$ respectively:
\begin{displaymath}
_{z}R^{k}=p(q^{-1}(z))\qquad
R^{k}_{y}=q(p^{-1}(y))
\end{displaymath}
$R^{k}_{y}$ is the set of points  on lines from  $ y$ 
intersecting  Y  with multiplicity
$\ge k$ and it is a cone 
  of vertex $y$. 
 In a neighborhood of $y$ we can identify  $\pro $ with  ${\bf C}^{m+r}$ where $y$ is the  origin,
$Y$ is defined  in an appropriate neighborhood of $y$ by $(r+1)$ polynomials
 $f_{1}\ldots f_{r+1}$. \newcommand{\ri}{R^{k}_{y}}             
 $\ri$\quad is given by vanishing of the homogeneous components of degree
$\le k-1$; and so if a generic line $ L $ of $\pro$ meets $\ri$  in $ k$ points, 
then $L\in \ri$. Moreover:
$$\textrm{dim} \ri \ge m+r-(r+1)(k-1) \forall y\in Y.$$
Let $F=q(R^{j+1})$ be the set of points of $\pro$ 
 on lines which intersect  $Y$ with multiplicity 
 $\ge j+1$ in one point: we want to prove that  $F$ is the hypersurface we looked for.
\\$Y\subset F$ because dim $R^{j+1}_{y}\ge 0$ $ \forall y\in Y.$
 The first  step is to prove that
$$F \subsetneqq  \pro.$$
Let $Y'=X\cdot G$ where $G$ is a generic hypersurface of degree $j$. 
$Y'$ is obtained by $Y$ by semicontinuity, so the dimension of $F$ 
passing from $Y$ to $Y'$ cannot decrease, and since the  (j+1)-secants to $Y'$ 
are contained in $G$  we obtain:
$\textrm{ dim} F\le m+r-1.$\\

Next step is to prove that :
$$\textrm{dim}F\ge m+r-1.$$

The set of  (j+1)-secants to $Y'=X\cdot G$ through an external point  
$P\in G$ is given by the intersection of  $\Sigma_{j+1}$ with the variety of 
Lemma 2.3, and so, by the assumption, we obtain a variety with degree different to $0$; $j$ times this degree
gives the virtual degree of $(j+1)$ secants intersecting a generic line of
 $\pro$. Since $Y$ is a degeneration of $Y'$, this virtual degree is the same 
and it is different from $0$ as stated previously.\\ 
Let $B$ the locus  of $(j+1)$-secants to $Y$ interecting a generic line, $B$
has dimension $\ge0$ in the  grassnammian of lines in ${\bf P}^{m+r}$ and it is  given by ${A\cap S}$ where:
$$A=\left\{\textrm{lines of $\pro$ that are $(j+1)$-secant to $Y$}\right\}$$
$$S=\left\{\textrm{lines of $\pro$ intersecting a given line} \right\}$$
$$\textrm{dim}\left\{A\cap S\}\ge 0 \Longrightarrow \textrm {codim} 
\{A\cap S\right\}      \le 2(m+r-1).$$
Since the line is generic, we have:
\[\textrm{codim}\{A\cup S\}=\textrm{codim A} +\textrm{ codim S}\]
\[\textrm{codim S} = m+r-2 \Longrightarrow \textrm{ codim}  A\le m+r\]
\[\textrm{dim} A\ge m+r-2\].
 Let $A'$ be the variety of points of $A$,
then  we have: $$\textrm{ dim} A'\ge m+r-1.$$
Now we have  to prove that $A'=F$.\\ 
The inclusion $F\subset A'$ is trivial; we want to prove that if  $p\in A'$, 
then $p \in F$. From Lemma $2.2$ we have that if  $p$ lies on a $(j+1)$-secant to $Y$ then 
it lies also on a line intersecting $Y$ with multiplicity  $(j+1)$ in a  point of $Y$.
Finally we have to prove that
$$deg F\le j.$$
Let suppose  that a generic line L of $\pro$ meets F in (j+1) points $z_{1}\ldots z_{j+1} \in L\cap F$.
Let's compute $c_{i}$= codim $(_{z_{i}}R^{j+1},Y)$:
$$\textrm{dim } R^{j+1} = \textrm{dim }Y + \textrm{dim }p^{-1}(y) = \textrm {dim} F + \textrm{dim }q^{-1}(z),$$         
since dim $ p^{-1}(y)=$ dim $R^{j+1}_{y}$ and dim $q^{-1}(z)=$dim $ _{z}R^{j+1}$,
 as  previously stated, we have:
$$\textrm{dim } R^{j+1}\ge 2m-1+r-(r+1)j$$
then
$$\textrm{dim } _{z_{i}}R^{j+1}\ge m-(r+1)j$$
$$c_{i}=\textrm{codim }(_{z_{i}}R^{j+1},Y)\le (r+1)j-1$$
By the Lefschetz-Barth's theorem  and by the assumption. we have
\begin{displaymath}
{\bf C}=  H^{2c_{i}}(\pro,{\bf C})=H^{2c_{i}}(Y,{\bf C})\Longrightarrow
  \bigcap_{i=1}^{j+1} {}_{z_{i}}R^{j+1}\neq\emptyset
\end{displaymath}

in fact: 
\begin{displaymath}
2c_{i}\le m-r-2 \quad e \quad(j+1)((r+1)j-1)\le m-1.
\end{displaymath}
Let $y\in\bigcap_{i=1}^{j+1} {}_{z_{i}}R^{j+1}$ 
then
$z_{i}\in L\cap R_{y}^{j+1} \quad per\quad i=1\ldots j+1$ 
and so 
$L\subset R_{y}^{j+1}.$
This is a contradiction as $L$ is generic.
We deduce that ${deg} F\le j$.

\section {Proof of theorem $1.2$}

{\bf Proof of theorem 1.2} 

Let $P$ be the fixed point and $Q\subset G({\bf P} ^{1},\pin )$ 
the space of lines from  $P$, $Q\simeq {\bf P} ^{n-1}$; let:
\begin{displaymath}
T=\{(q,l) \mid q\in {\bf P}^ {n}\quad l\in Q \quad q\in l\}
\end{displaymath}
and $\alpha$ and $\beta$ be the projections of  $T$ on $\pin$ and  $Q$ respectively.
\\$T$ is a ${\bf P}^{1}$-bundle on $Q$ and the fibre is given by all the points 
lying  on lines $l$, we can view  it as the projectivised of $\OO_{Q}\oplus \OO_{Q}(-1)$.
\\Let $(T/Q)^{k+1}$ be the $(k+1)$-power of fibre of $T$ on $Q$,that is:
\begin{displaymath}
(T/Q)^{k+1}=\underbrace{T\times_{Q}T\times_{Q}\ldots\times_{Q}T}_{ k+1  {\textrm times}}  
\end{displaymath}
We call $Z\in (T/Q)^{k+2}$ the incidence variety in 
$T\times_{Q}(T/Q)^{k+1}$, that is:
\begin{displaymath}
Z=\{(x_{0},\ldots,x_{k+1}\in (T/Q)^{k+2} \quad \mid x_{0}=x_{i} \textrm{ for same } i\in (1, \ldots, k+1)\}
\end{displaymath}
Let $p$ and $q$ be the projections of $Z$ on  $T$ and $(T/Q)^{k+1}$ respectively; 
we denote:
\begin{displaymath}
E^{(k+1)}=q_{*}(p^{*}\alpha^{*}(E))
\end{displaymath}
$E^{(k+1)}$ is a vector bundle on $(T/Q)^{k+1}$
of rank $r(k+1)$.
 Let $s$ be the section of $E$ such that $X$ is the zero locus of $s$; 
$s^{(k+1)}=q_{*}(p^{*}\alpha^{*}(s))$
is a section of  $E^{(k+1)}$ which vanishes in the set:
$\{(x_{1},\ldots, x_{k+1}\in (T/Q)^{k+1} \quad \mid \alpha(x_{i})\in X \}$.
The line through $\alpha(x_{1}),\ldots, \alpha(x_{k+1})$ is a $(k+1)$-secant to $X$.
Considering that the rearrangement of those points  gives the same
$(k+1)$-secant to $X$, from  Portous' formula we have that, if the dimension  
is  zero, the number of $(k+1)$-secant is given by the degree of the top 
Chern-class $c_{(k+1)r}(E^{(k+1)})$ divided by $(k+1)!$.
So if we want to know the degree of $(k+1)$-secants we have to compute  $c_{(k+1)}(E^{(k+1)})$.
Let $q_{1},\ldots q_{k+1}$ be the  projections of $(T/Q)^{k+1}$ on $\pin$; 
we have the following exact sequence:
\begin{displaymath} 
0\longrightarrow q^{*}E\otimes \OO(-\Delta_{1,k+1}\ldots  -\Delta_{k,k+1})\longrightarrow E^{(k+1)} \longrightarrow E^{(k)} \longrightarrow 0
\end{displaymath}
with
$\Delta_{i,j}=\{(x_{1},\ldots x_{k+1}\in (T/Q)^{k+1} \quad \mid x_{i}=x_{j}\}$.
From sequence we have:
\begin{displaymath}
c_{(k+1)r}E^{(k+1)}=c_{kr}E^{(k)} c_{r}(q^{*}E\otimes \OO(-\Delta_{1,k+1}\ldots -\Delta_{k,k+1})).
\end{displaymath}
It is necessary to determine the  cohomology of $T$ and $(T/Q)^{k+1}$.\\
Let  $\alpha$ and $\beta$ be the projections of $T$  on $\pin$ and  $Q$ respectively: 
$T$ is blow-up of $\pin$ in $P$; we call $D$ 
the exceptional divisor and $H= \alpha^{*}(\OO_{\pin}(1))$, 
then we have:
$H-D=\beta^{*}(\OO_{Q}(1)), D=\OO_{T}(1).$
The  Wu-Chern's equation gives:
$D^{2}+\beta^{*}\OO_{Q}(1)D=0.$\\
The intersection ring of  $T$ is generated by two elements:
\begin{displaymath}
\langle D,H-D \rangle =\langle D, \beta^{*}\OO_{Q}(1)\rangle.
\end{displaymath}
\newcommand{\oq}{\beta^{*}\OO_{Q}(1)}
For the next degrees we have:
\begin{displaymath}
(\oq)^{2}D=-\oq D^{2}=D^{3}
\end{displaymath}
$$\vdots $$
\begin{displaymath}
 (\oq)^{n}D=(-1)^{n-1}\oq D^{n}=D^{n+1}.
\end{displaymath}
 We observe that 
$H^{n}=1$ and $ D^{n}=(-1)^{n-1}$
in fact $D_{\mid D}=\OO_{D}(-1)$ and $D\simeq {\bf P}^{n-1}$.
Consider now the fibred product $T\times_{Q} T$:
$H^{*}(T)$ is  generated by $D$ as $H^{*}(Q)$-module;
$H^{*}(T\times_{Q}T)=H^{*}(T)\times_{H^{*}(Q)} H^{*}(T)$ is 
generated by $D\otimes 1=D_{1}$ and $1\otimes D= D_{2}$ as 
$H^{*}(Q)$-module; if we consider it as a vector space we have:
\begin{displaymath}
H^{2}(T\times_{Q}T)=\langle D_{1},D_{2},\oq \rangle
\end{displaymath}
\begin{displaymath}
H^{4}(T\times_{Q}T)=\langle D_{1}^{2},D_{2}^{2},(\oq) ^{2},D_{1}D_{2}\rangle
\end{displaymath}
$$\vdots $$
\begin{displaymath}
H^{2j}(T\times_{Q}T)=\langle D_{1}^{j},D_{2}^{j},(\oq) ^{j},D_{1}D_{2}\oq ^{j-2} \rangle
\end{displaymath} 
 we denote:
$H_{1}=q_{1}^{*}(H) \quad H_{2}=q_{2}^{*}(H)$
\[H_{1}=D_{1}+\oq \quad H_{2}=D_{2}+\oq.\]
We prove that:
$\Delta_{1,2}=D_{1}+D_{2}+\oq$ .
Let $p_{1}$ and $p_{2}$ be the projection of $T\times_{Q}T$ on the two factors.
$\Delta_{1,2}$ is given by the zero locus of a section of $p^{*}_{1}\OO(1)\otimes p^{*}Q_{rel}$ (see \cite{O-S-S},page 242).
We know that $c_{1}(p^{*}_{1}\OO(1))=D_{1}$; consider now the exact section
$$ 0\longrightarrow \OO (-1)\longrightarrow \beta^{*}\OO \oplus \beta^{*}\OO (1)\longrightarrow Q_{rel}\longrightarrow 0.$$
We get $c_{1}(p^{*}_{2}Q_{rel})=D_{2}+\beta^{*}\OO(1)$ and so we have $\Delta_{1,2}=D_{1}+D_{2}+\oq .$ 

For the general case we have that:\\
$H^{2}(T/Q)^{k+1}$ is generated by  $D_{1},D_{2},\ldots D_{k+1},\oq $,\\
\vdots\\
$H^{2m}(T/Q^{k+1})$ is generated by $(\oq )^{m},D_{i_{1}}\ldots D_{i_{t}}\oq ^{m-t}.$
Moreover:
$\Delta_{i,j}=D_{i}+D_{j}+\oq.$

Now we prove the theorem proceeding by induction on $k$:
for $k=1$ the exact sequence is:
 $$0\longrightarrow q_{2}^{*}(E) \otimes O(-\Delta_{1,2})\longrightarrow E^{(2)}\longrightarrow q_{1}^{*}(E)\longrightarrow 0$$

\begin{eqnarray*}
c_{2r}(E)^{(2)}&=&c_{r}((q_{1}(E))c_{r}(q_{2}^{*}(E)\otimes  O(-D_{1}-D_{2}-\oq)\\
 &=&c_{r}(E)H_{1}^{r}[c_{r}(E)H_{2}^{r}+c_{r-1}(E)H_{2}^{r-1}(-D_{1}-D_{2}-\oq)+\\
& &\ldots + c_{r-i}(E)H_{2}^{r-i}(-D_{1}-D_{2}-\oq)^{i}+\\
& &\ldots +(-D_{1}-D_{2}-\oq)^{r}]
\end{eqnarray*}
since: $H_{1}D_{1}=0$ and $D_{2}+\oq =H_{2}$
we have:
$c_{2r}=c_{r}(E)c_{r}(-1)H_{1}^{r}H_{2}^{r}$. 
Now we suppose  the statement true for $n\ge k$ 
and we try to prove it for  $n=k+1$.
\begin{eqnarray*}
c_{(k+1)r}&=&c_{kr}E^{(k)}c_{r}(q_{k+1}^{*}(E)\otimes O(-D_{1}-D_{2}\ldots -D_{k}-kD_{k+1}-k\oq)\\
&= &c_{r}(E)c_{r}(E(-1))\ldots c_{r}(E(-k+1))H^{r}_{1}\ldots H_{k-1}^{r}[c_{r}(E)H_{k+1}^{r}+\\
& & c_{r-1}(E)H_{k+1}^{r-1}(-D_{1}-D_{2}\ldots)+\ldots (-D_{1}\ldots-k\oq)^{r}]
\end{eqnarray*}
since we know that: 
$H_{i}D_{i}=0$
and $D_{k+1}+\oq =H_{k+1}$
we obtain:
\begin{displaymath}
c_{(k+1)r}(E^{(k+1)})=c_{r}(E)c_{r}(E(-1))\ldots c_{r}(E(-k))H_{1}^{r}\ldots H_{k+1}^{r}
\end{displaymath}
and so the theorem is proved. 
\\
{\bf Remark 1} We observe that  if the dimension of the locus of $k$-secants 
through a generic point is smaller than expected, then the class of the formula 
has to be zero (see \cite{G} Remark $2.2$).\\

{\bf Remark 2} In the case $r=2$ the theorem 
has been already proved by 
Ran in [R]. By the Hartshorne-Serre correspondence  every subcanonical subvariety of codimension $2$   is a zero locus of a section of a rank $2$ vector bundlon ${\bf P}^{n}$; moreover if $n\ge 10$ by Larsen's theorem we have that every subvariety is subcanonical. In this case the  formula for  $j+1$-secant is true
for every subvariety.\\
\section{A new proof of Zak theorem on linear normality}
Let $X$ be a  $r$ codimensional subvariety  of $\pin$; from  Barth theorem 
we have that if $r\ge n/4$ then $H^{2i}(X,\ZZ)\simeq \ZZ$; in particular we can write  $c_{i}(N)= c_{i}H^{i}$ with $i=1 
\ldots r$ where $c_{i}\in {\bf Z}$.
From now, we consider  $c_{i}(N)$ as a integer.
\begin{lemma}
Let $X$ be  a  $r$ codimensional subvariety  of $\pin$. If $n\ge 4r$, then the degree of set of 
bisecant to $X$  through an external point  is 
$$c_{r}(N)c_{r}(N(-1)).$$
\end{lemma}
{\bf Proof.} Let $P$ be the fixed point, if we project $X$ from $P$ to a generic hyperplane  
we can use the  double point formula \cite{F}  to get the set of bisecant to $X$ from $P$. 
\[2\Sigma_{2}=f^{*}f_{*}[X]-(c(f^{*}T{\bf P}^{n-1})c(TX)^{-1})_{r-1}\cap  [X]\]
and from the exact sequence
\[0\longrightarrow T_{X}\longrightarrow T{\bf P}^{n}_{\mid S}\longrightarrow N_{X,{\bf P}^{n}}\longrightarrow 0\]    
we have $ c(TX)^{-1}=c(T{\bf P}^{n})^{-1} c(N)$ and substituting we get
$$2\Sigma_{2}=H^{r-1}(d-c_{r-1}+c_{r-2}+\dots (-1)^{i}c_{i}\dots)= c_{r}(N(-1))H^{r-1}. $$
From theorem $1.1$ and from Lemma we get a different proof of Zak's theorem.
\begin{thm}[Zak]
Let $X$ be a $r$ codimensional subvariety of $\pin$, if $n\ge 4r$, then $X$ is linearly normal.
\end{thm}
{\bf Proof.}
We prove the theorem proceeding by induction on $r$. If $r=1$ is trivially true. Now we suppose that it is true for $r-1$. If $c_{r}(N(-1))\neq 0$ from theorem $1.1$ and from lemma $4.1$ we have the thesis. If $c_{r}(N(-1))=0$ from lemma $4.1$ we have that there are not bisecant to $X$ through an external point $P$; if we project $X$ from $P$ to a generic hyperplane, we get a smooth subvariety in ${\bf P}^{n-1}$ of codimension $r-1$ that is linearly normal by induction.
This is a contradiction. 

\section{Proof of theorem $1.4$}
\begin{lemma}
Let $l,m,p \in {\bf N}$ such that $l+p=m$ then we have 
$$\binom{l}{t}=\sum_{i=0}^{k}(-1)^{i}\binom{m}{t-i}\binom{p-1+i}{i}$$
 $$\sum_{1=0}^{t}(-1)^{i}\binom{n}{ t-i}\binom{n+1+i}{ i}=(-1)^{t}$$ 

\end{lemma}
{\bf Proof.} We consider an exact sequence 
$$0\longrightarrow A \longrightarrow B\longrightarrow C \longrightarrow 0 $$         
where  A, B, C  are vector spaces with dimension respectively  $l,m,p$. 
From this exact sequence we obtain other two exact sequences:
\begin{equation} 
0\longrightarrow \wedge ^{t} A \longrightarrow \wedge^{t} B\longrightarrow 
\wedge^{t-1}B\otimes C\dots \longrightarrow  \wedge^{t-i}\otimes S^{i}C 
\longrightarrow \dots \rightarrow S^{n}C \longrightarrow 0 
\end{equation}

\begin{equation}                       
0 \longrightarrow \wedge^{t} A\longrightarrow \wedge^{t-1}A\otimes B\dots 
\longrightarrow  \wedge^{t-i}A\otimes S^{i}B \rightarrow \dots \rightarrow S^{n}B
\rightarrow S^{n}C\rightarrow  0 
\end{equation}
considering that $\wedge^{t}({\bf C}^{m})=\binom{m}{ t} $ and 
$S^{t}({\bf C^{n}})=\binom{n-1+t}{ t}$
we have
$$\binom{l}{ t}=\sum_{1=0}^{t}(-1)^{i}\binom{m}{ t-i}\binom{p-1+i}{ i}$$
From $(2)$  if $m=n+1$ and $l=n$ we have:
 $$\sum_{1=0}^{t}(-1)^{i}\binom{n}{ t-i}\binom{n+1+i}{ i}=(-1)^{t}$$ 

\begin{lemma}
Let $X$ a  $r$ codimensional subvariety  of $\pin$ then if $n\ge 4r$ the  locus of trisecant  is 
$$\Sigma_{3}=\frac{1}{2}H^{2r-2}c_{r}(N(-1))c_{r}(N(-2))$$
\end{lemma}

{\bf Proof}
 G\"{o}ttsche's formula for trisecant through a fixed point  is 
$$\Sigma_{3}=(a)+(b)-(c)$$
where:
$$(a)=H^{2r-2} \left(\frac{n}{2}d^{2}-\sum_{k=0}^{n-r}\left(\binom{2n-2r+2}{ k}-
\binom{n}{ k-n+2r-2}\right)\int_{X}H^{k}s_{n-r-k}/2\right)$$ 
$$(b)=\sum_{k=0}^{2r-2}\sum_{t=0}^{n-1}\binom{n}{ t}\binom{n+1}{ k-t}
\sum_{j=r-t-1}^{2r-2-k}2^{j+t-r+1}s_{j}(X)s_{2r-2-k-j}(X)H^{k}$$
 and
$$(c)= \sum_{k=0}^{2r-2}d\binom{n+r}{ k}s_{2r-2-k}(X)H^{k}.$$

We prove the Lemma when $r$ is even  (the case $r$ odd is the same).
It is well known:
$$s_{k}=\sum_{i=0}^{n}(-1)^{k+1}H^{k-i}c_{i}(N)\binom{n+k-i}{ k-i}.$$
Let $c_{i}=c_{i}(N)$; substituting we have:
\begin{eqnarray*}
(a)&=&H^{2r-2} \Big(\frac{n}{2}d^{2}-\frac{1}{2}
\Big(\sum_{k=0}^{n-r}\binom{2n-2r+2}{ k}\cdot\\
& &\sum_{i=0}^{r}(-1)^{n-r-k+i}c_{i}H^{n-r-i}
\binom{n+n-r-k-i}{ n-r-k-i}+\\
& & -\sum_{k=0}^{n-r}\sum_{i=0}^{r}(-1)^{n-r-k+i}\binom{n}{ k-n+2r-2}
\binom{n+n-r-k-i}{ n-r-k-i}\Big)\Big)
\end{eqnarray*}
we put $k'=n-r-k-i$ and so we have:
\begin{eqnarray*}
(a)&=&H^{2r-2} \Big(\frac{n}{2}d^{2}-\frac{1}{2}
\Big(\sum_{i=0}^{r} c_{i}H^{n-r-i}
\sum_{k'=0}^{n-r-i}(-1)^{k'}\binom{n+k'}{ k'} \binom{2n-2r+2}{ n-r-i-k'}+\\
& &-\sum_{k=0}^{r-2-i}(-1)^{k'}\binom{n}{ r-2-i-k'}\binom{n+k'}{ k'}\Big)\Big)
\end{eqnarray*}
now we can use the Lemma $5.1$ and we obtain
\begin{equation}(a)=H^{2r-2} \left(\frac{n}{2}d^{2}-\frac{1}{2}\left(d^{2}(n-2r+1)-\sum_{i=0}^{r-1}(-1)^{r-i}c_{i}H^{n-r-i}\right) \right)
\end{equation}
\begin{eqnarray*}
(b)&=&\sum_{t=0}^{n-1}\binom{n}{ t}\sum_{j=r-t-1}^{2r-2} 2^{j+t+1-r}s_{j}
\sum_{k=t}^{2r-2-j}\binom{n+1}{ k-t}\sum_{m=0}^{r} (-1)^{2r-2-k-j+m}\cdot\\
& & H^{2r-2-j-m}
 c_{m}\binom{n+2r-2-k-j-m}{ 2r-2-k-j-m}
\end{eqnarray*}
If we denote  $k'=2r+2-k-j-m$ we get 
\begin{eqnarray*}
(b)&=&\sum_{m=0}^{r}c_{m}\sum_{t=0}^{n-1}\binom{n}{ t}\sum_{j=r-t-1}^{2r-2} 2^{j+t+1-r}s_{j}
H^{2r-2-j-m}\cdot \\
& & \sum_{k'=0}^{2r-2-j-t-m}(-1)^{k'}
\binom{n+1}{2r-2-j-m-t-k'} 
\binom{n+k'}{ k'}
\end{eqnarray*}
From  Lemma $5.1$ we have that  the last sum is equal to $1$ if $2r-2-j-m-t=0$ and 
equal to $0$ in the other cases; this fact implies also that $j=2r-2-m-t\ge r-t-1$ and so we obtain $m\le r-1$.

$$(b)=\sum_{m=0}^{r-1}c_{m}\sum_{t=0}^{n-1}\binom{n}{t} 2^{r-1-m}s_{2r-2-m-t}
H^{t}$$
\begin{eqnarray*}
(b)&=&\sum_{m=0}^{r-1}c_{m}2^{r-1-m}\sum_{t=0}^{n-1}\binom{n}{ t} 
\sum_{i=0}^{r} (-1)^{2r-2-m-t+i}\cdot\\
& &H^{2r-2-i-m}c_{i}\binom{n+2r-2-m-t-i}{ 2r-2-m-t-i}
\end{eqnarray*}
let $t'=2r-s-m-t-i$
\begin{eqnarray*}
(b)&=&\sum_{m=0}^{r-1}\sum_{i=0}^{2r-2-m}c_{m}c_{i}H^{2r-2-i-m}2^{r-1-m}\cdot\\
& &\sum_{t'=0}^{2r-2-m-i}
 (-1)^{t'}\binom{n+t'}{ t'}\binom{n}{ 2r-2-m-t-i}
\end{eqnarray*}
and again from the Lemma $5.1$ we have
\begin{equation}
(b)=\sum_{m=0}^{r-1}\sum_{i=0}^{2r-2-m}(-1)^{m+i}2^{r-1-m}c_{m}c_{i}H^{2r-2-i-m}
\end{equation}
$$(c)= \sum_{k=0}^{2r-2}d\binom{n+r}{ k}\sum_{i=0}^{r} (-1)^{2r-2-k+i}H^{2r-2-i}c_{i}\binom{n+2r-2-k-i}{ 2r-2-k-i}$$
let $k'=2r-2-k-i$ 
\begin{eqnarray*}
(c)&=&\sum_{i=0}^{r} d H^{2r-2-i}c_{i} \sum_{k'=0}^{2r-2}\binom{n+r}{ 2r-2-i-k'} (-1)^{k'}
\binom{n+k'}{ k'}\\
&=&\sum_{i=0}^{r} d H^{2r-2-i}c_{i}  \binom{r-1}{ 2r-2-i}
\end{eqnarray*}
\begin{equation}
(c)= d c_{r-1}  H^{r-1}+d^{2}(r-1) H^{2r-2}.
\end{equation}
Supposing that we are in the range of Barth's theorem, we have that $c_{i}=c_{i}H^{i}$
where $ c_{i}\in \ZZ $.
Finally we get from (3), (4) and (5):
$$\Sigma_{3}=(a)+(b)-(c)=\frac{1}{2}H^{2r-2}\sum_{m=0}^{r}\sum_{i=0}^{r}(-1)^{m+i}2^{r-m}c_{m}c_{i}$$ 
that is 
$$\Sigma_{3}=H^{2r-2} \frac{1}{2}c_{r}(N(-1))c_{r}(N(-2))$$.

\end{document}